\documentclass[a4paper, 10pt,oneside]{amsart}
\usepackage{amsmath, amssymb, amsfonts}
\usepackage[all]{xy}
\newtheorem{tm}{Theorem}
\newtheorem{lm}[tm]{Lemma}
\newtheorem{prop}[tm]{Proposition}

\newtheorem{napomena}{Remark}
\newenvironment{dokaz}
{\noindent\emph{Proof:}\ }
{\hfill $\blacksquare$}

\newenvironment{sketch}
{\noindent\emph{Sketch of proof:}\ }
{\hfill $\blacksquare$}

\newcommand{\Z}
{{\mathbb Z}}
\newcommand{\N}
{{\mathbb N}}
\newcommand{\C}
{{\mathbb C}}

\newcommand{\g}
{{\mathfrak g}}

\newcommand{\gt}
{\tilde{{\mathfrak g}}}

\newcommand{\h}
{{\mathfrak h}}
\newcommand{\n}
{{\mathfrak n}}
\newcommand{\gsl}
{{\mathfrak sl}}

\newcommand{\Gamt}
{\tilde{{\Gamma}}}
\newcommand{\Weyl}
{{\mathcal W}}      

\begin{document}

\title[Bases of Feigin-Stoyanovsky's type subspaces
for $C_\ell^{(1)}$]{Bases of Feigin-Stoyanovsky's type subspaces
for  $C_\ell^{(1)}$}
\author{Ivana Baranovi\' c, Mirko Primc and Goran Trup\v cevi\' c}

\email{baranovic.ivana@yahoo.com, primc@math.hr, goran.trupcevic@ufzg.hr}

\address{Ivana Baranovi\' c, University of Zagreb, Faculty of Chemical Engineering and Technology, Zagreb, Croatia}

\address{Mirko Primc, University of Zagreb, Faculty of Science, Department of Mathematics, Zagreb, Croatia}

\address{Goran Trup\v cevi\' c, University of Zagreb, Faculty of Teacher Education,  Zagreb, Croatia}

\subjclass[2000]{Primary 17B67; Secondary 17B69, 05A19.\\
\indent Partially supported by Croatian Science Foundation under the project 2634 and by the Croatian Scientific Centre of Excellence QuantixLie}

\begin{abstract}
In this paper we construct combinatorial bases of Feigin-Stoyanovsky's type subspaces of standard modules
for level $k$ affine Lie algebra $C_\ell^{(1)}$.
We prove spanning by using annihilating field $x_\theta (z)^{k+1}$ of standard modules. 
In the proof of linear independence we use simple currents and intertwinining operators whose existence
is given by fusion rules.
\end{abstract}
\maketitle

\section{Introduction}

Let $\g$ be a simple complex Lie algebra, $\h\subset\g$ its Cartan
subalgebra, $R$ the corresponding root system. Let $\g=\h+\sum_{\alpha\in R}\g_\alpha$ be a root
decomposition of $\g$. Fix root vectors $x_\alpha\in\g_\alpha$. Let
\begin{equation}\label{E: ZGradG}
\g=\g_{-1}\oplus \g_0 \oplus \g_1
\end{equation} 
be a $\Z$-gradation of $\g$, where $\h\subset\g_0$. Such gradations correspond to a choice of a 
minuscule coweight $\omega\in\h$. Denote by $\Gamma\subset R$ a set
of roots such that $\g_1=\sum_{\alpha\in\Gamma}\g_\alpha=\sum_{\omega(\alpha)=1}\g_\alpha$. We call $\Gamma$ the {\em set of colors}.

Affine Lie algebra associated with $\g$ is $\gt=\g\otimes
\C[t,t^{-1}]\oplus \C c \oplus \C d$, where $c$ is the canonical
central element, and $d$ the degree operator. Elements
$x_\alpha(n)=x_\alpha\otimes t^n$ are fixed real root vectors.
The $\Z$-gradation of $\g$ induces analogous gradation of $\gt$:
$$\gt=\gt_{-1}\oplus \gt_0 \oplus \gt_1,$$
where $\gt_1=\g_1\otimes\C[t,t^{-1}]$ is a commutative Lie
subalgebra with a basis
$$\Gamt=\{x_\gamma(n)\,|\,n\in\Z,\gamma\in\Gamma\}.$$

For a standard $\gt$-module $L(\Lambda)$ of level $k=\Lambda(c)$,
define a {\em Feigin-Stojanovsky's
type subspace} $W(\Lambda)$  as a $\gt_1$-submodule generated with
a highest weight vector $v_\Lambda$,
$$W(\Lambda)=U(\gt_1)\cdot v_\Lambda\subset L(\Lambda).$$

In this paper, for a Lie algebra $\g$ of type $C_\ell$, we construct a basis of $W(\Lambda)$ consisting of monomial vectors $x(\pi)v_\Lambda$, 
where $x(\pi)$ are monomials in $\Gamt$. 
Poincare-Birkhoff-Witt's theorem gives a spanning set of  $W(\Lambda)$ 
\begin{equation} \label{E: MonVect}
\{x_{\gamma_1}(-n_1)
x_{\gamma_2}(-n_2)\cdots x_{\gamma_t}(-n_t) v_\Lambda\,|\,
t\in\Z_+,\gamma_i\in \Gamma, n_i\in\N\}.
\end{equation}
In order to obtain a basis of $W(\Lambda)$ we find relations for standard modules upon which we reduce the spanning set.
Finally, we prove linear independence by	using intertwining operators.


The notion of Feigin-Stoyanovsky's type subspaces is similar to a notion of {\em principal subspaces} that were introduced by 
B.\,Feigin and A.\,Stoyanovsky  for $\gsl_2(\C)$ and $\gsl_3(\C)$ (\cite{FS}). In this case one looks at a triangular decomposition of $\g$
instead of \eqref{E: ZGradG}. Many different authors have studied these spaces,
their bases, character formulas, exact sequences etc. (\cite{Ge}, \cite{Sto}, \cite{CLM1}, \cite{CLM2}, \cite{AKS}, \cite{C}, \cite{CalLM1}, \cite{CalLM2}, \cite{Bu}, \cite{Sa}).


Another type of principal subspaces, the so called Feigin-Stoyanovsky's type subspaces $W(\Lambda)$ defined above, was implicitly studied in \cite{P1}
for $\gsl_{\ell+1}(\C)$. It turned out that in this case, bases are parameterized by $(k,\ell+1)$-admissible configurations, studied by Feigin et al.\,(\cite{FJLMM}, \cite{FJMMT}).
The $\Z$-gradations \eqref{E: ZGradG} are closely related to simple current operators (\cite{DLM}). We hope that this kind of construction of combinatorial bases 
will be possible for all affine Lie algebras. Up to now, this was done for the type $A_\ell^{(1)}$ in \cite{P3}, \cite{T1} and \cite{T2}, for $B_2^{(1)}$ in \cite{P4},
for $D_4^{(1)}$, level $1$ and $2$, in \cite{Ba}, and for all classical types, level $1$, in \cite{P2}.

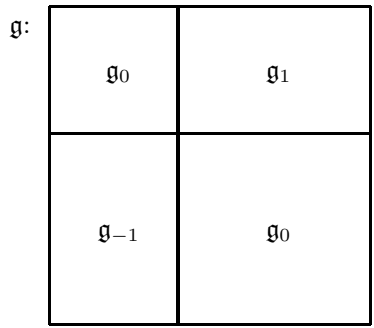
\begin{figure}[ht] \caption{$\Z$-gradation of $\gsl_{\ell+1}(\C)$}\label{ZGrad_fig}
	\begin{center}\begin{picture}(145,140)(-15,-10)
		\thicklines \put(0,0){\line(0,1){120}} \put(0,0){\line(1,0){120}}
		\put(48,0){\line(0,1){120}} \put(0,72){\line(1,0){120}}
		\put(120,120){\line(0,-1){120}} \put(120,120){\line(-1,0){120}}
		\put(-15,110){$\g$:} \put(21,93){$\g_0$} \put(81,33){$\g_0$}
		\put(18,33){$\g_{-1}$} \put(81,93){$\g_1$}  \linethickness{.075mm}
		\end{picture}\end{center}
\end{figure}

 The first general case beyond admissible cofigurations was given in \cite{T1} and \cite{T2} 
 where a new kind of combinatorial conditions emerged.
 The minuscule coweight $\omega$ that was considered in these papers corresponds to a $\alpha_m$, $1\leq m \leq \ell$.
If $\g=\gsl_{\ell+1}(\C)$ represents matrices of trace $0$ then in the $\Z$-gradation \eqref{E: ZGradG} the subalgebra ${\mathfrak g}_0$ 
consists of block-diagonal matrices, while ${\mathfrak g}_1$ and $ {\mathfrak g}_{-1}$ consist of
matrices with non-zero entries only in the upper right or lower-left block, respectively. This is illustrated on the figure \ref{ZGrad_fig}.

The set of colors $\Gamma$ corresponds to a rectangle with rows $1,\dots,m$ and columns $m,\dots,\ell$.
Monomial basis of $W(\Lambda)$ is described in terms of {\em difference} and {\em initial conditions}.
 
In the level $1$ case a monomial vector \eqref{E: MonVect} satisfies difference conditions for $W(\Lambda_r)$, if 
colors of elements of the same degree $-n$ lie on a diagonal path in $\Gamma$ as shown in figure \ref{DCSeqA_fig}.
Furthermore, colors of elements of  degree $-n-1$ lie on a diagonal path that lies below or to the left of the 
preceding path.

\begin{figure}[ht] \caption{Difference conditions for $\gsl_{\ell+1}(\C)$ level $1$} \label{DCSeqA_fig}
	\begin{center}\begin{picture}(200,140)(-8,-10) \thicklines
		\put(0,0){\line(1,0){180}} \put(0,0){\line(0,1){120}}
		\put(180,120){\line(0,-1){120}} \put(180,120){\line(-1,0){180}}
		\put(-6,111){$\scriptstyle 1$} \put(-6,99){$\scriptstyle 2$}
		\put(-8,3){$\scriptstyle m$} \put(3,-8){$\scriptstyle m$}
		\put(174,-8){$\scriptstyle \ell$}
		
		\linethickness{.075mm} \multiput(98,77)(4,0){21}{\line(1,0){2}}
		\multiput(97,78)(0,4){11}{\line(0,1){2}}
		
		\thinlines \put(96,76){$\scriptscriptstyle \bullet$}
		\put(98,78){\line(1,1){12}} \put(108,88){$\scriptscriptstyle
			\bullet$} \put(132,94){$\scriptscriptstyle \bullet$}
		\put(134,96){\line(1,1){12}} \put(144,106){$\scriptscriptstyle
			\bullet$} \put(110,90){\line(4,1){24}} \put(108,97){$\scriptstyle
			(-n)$}
		
		\put(60,10){$\scriptscriptstyle \circ$}
		\put(84,16){$\scriptscriptstyle \circ$}
		\put(108,40){$\scriptscriptstyle \circ$}
		\put(120,64){$\scriptscriptstyle \circ$}
		\put(62,12){\line(4,1){24}} \put(86,18){\line(1,1){24}}
		\put(110,42){\line(1,2){12}} \put(96,22){$\scriptstyle (-n-1)$}
		\end{picture}\end{center}
\end{figure}

A monomial vector \eqref{E: MonVect} satisfies initial conditions for $W(\Lambda_r)$ if 
a diagonal path of colors of degree $-1$ lie either below the
$r$-th row (if $1\leq r \leq m$), or on the left of the $r$-th column (if $m\leq r \leq \ell$).

In the case considered in this paper, when $\g$ is a simple Lie algebra  of type $C_\ell$,  similar combinatorics appear.
The set of colors $\Gamma$ can be represented as a triangle with rows and columns ranging from $1$ to $\ell$ (cf. figure \ref{Gamma_fig}).

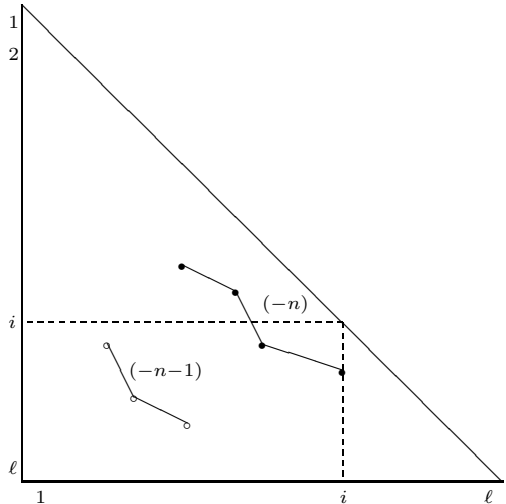
\begin{figure}[ht] \caption{Difference conditions for $C_\ell$ level $1$} \label{DClev1_fig}
	\begin{center}\begin{picture}(200,200)(-10,-10) 
		\put(0,0){\line(1,0){180}} \put(0,0){\line(0,1){180}}
		\put(180,0){\line(-1,1){180}} 
		\put(-5,171){$\scriptstyle 1$}
		\put(-5,159){$\scriptstyle 2$}
		\put(-5,3){$\scriptstyle \ell$} \put(5,-8){$\scriptstyle
			1$}
		\put(173,-8){$\scriptstyle \ell$}
		\put(90,65){$\scriptstyle (-n)$}
		\put(58,80){$\scriptscriptstyle \bullet$}
		\put(78,70){$\scriptscriptstyle \bullet$}
		\put(88,50){$\scriptscriptstyle \bullet$}
		\put(118,40){$\scriptscriptstyle  \bullet$}
		\put(60,82){\line(2,-1){20}}
		\put(80,72){\line(1,-2){10}}
		\put(90,52){\line(3,-1){30}} 
		\multiput(120,60)(-4,0){30}{\line(-1,0){2}}
		\multiput(120,60)(0,-4){15}{\line(0,-1){2}}
		\put(-5,58){$\scriptstyle i$}
		\put(119,-8){$\scriptstyle i$}
		\put(30,50){$\scriptscriptstyle \circ$}
		\put(40,30){$\scriptscriptstyle \circ$}
		\put(60,20){$\scriptscriptstyle \circ$}
		\put(32,52){\line(1,-2){10}} 
		\put(42,32){\line(2,-1){20}}
		\put(40,40){$\scriptstyle
			(-n-1)$}
		\end{picture}\end{center}
\end{figure}

For a fundamental weight $\Lambda_r$, a monomial vector \eqref{E: MonVect} satisfies difference conditions for $W(\Lambda_r)$, if 
colors of elements of two consecutive degrees lie on diagonal paths that are related in a way shown in figure \ref{DClev1_fig}, i.e. if $(-n)$-path ends in the $i$-th column, then $(-n-1)$-path lies below the $i$-th row.
A monomial vector \eqref{E: MonVect} satisfies initial conditions for $W(\Lambda_r)$ if 
a diagonal path of colors of degree $-1$ lies below the
$r$-th row.

In the proof of linear independence we follow ideas of Georgiev (cf. \cite{Ge}), and of Caparelli, Lepowsky and Milas (cf. \cite{CLM1}, \cite{CLM2}). 
Start from a relation of linear dependence
\begin{equation} \label{E: LinDep_Intro}
\sum c_{\pi}x(\pi) v_\Lambda = 0
\end{equation}
where $x(\pi)$ are some monomials that satisfy difference and initial conditions for $W(\Lambda)$ and $c_\pi\in\C$. 
The main idea 
is to use intertwinining operators between standard modules (cf. \cite{FHL}, \cite{DL}) and simple current operators (cf. \cite{DLM}) to reduce this relation to a relation of linear dependence on another 
standard module and proceed inductively. 

More concretely, let  $x(\mu)$ be, in some sense, the smallest monomial in \eqref{E: LinDep_Intro}. Then there is a 
coefficient of intertwining operator $w$ that commutes with $\gt_1$ and which sends $v_\Lambda$ to a vector $v'$ 
from the top of a standard module $L(\Lambda')$ that is annihilated by almost all monomials greater than $x(\mu)$. 
Furthermore, for the remaining monomials the action of $x(\pi)$ on $v'$ yields $x(\pi_2) [\omega] v_{\Lambda''}$, 
where $[\omega]$ is a simple current operator and $x(\pi_2)$ is a submonomial of $x(\pi)$. On the other side, commutation of a monomial
with $[\omega]$ raises degrees of factors by $1$; thus we obtain 
\begin{equation} \label{E: LinDep4_Intro}
0= \sum c_{\pi} [\omega] x(\pi_2^+) v_{\Lambda''},
\end{equation}
where $x(\pi_2^+)$ is obtained from $x(\pi_2)$ by raising degrees of factors by $1$.
Finally, simple current operator $[\omega]$ is a linear injection, hence
\begin{equation} \label{E: LinDep5_Intro}
0= \sum c_{\pi}  x(\pi^+) v_{\Lambda''}.
\end{equation}
This is a relation of linear dependence on $L(\Lambda'')$ with monomials of higher degree than in \eqref{E: LinDep_Intro}.
Since they are obtained from the ones from \eqref{E: LinDep_Intro} by raising their  degrees, it turns out that
they also satisfy difference and initial conditions, this time for $W(\Lambda'')$. This enables us to use inductive argument to obtain $c_\mu=0$. 

In the higher level case, difference and initial conditions are again given in terms of certain paths in $\Gamma$. Moreover, just like
in the $A_\ell$ case considered in \cite{T2}, in the $C_\ell$ case one can embed $L(\Lambda)$ into a tensor product of level $1$ modules and 
factorize monomial vectors from the basis into a tensor product of level $1$ monomial vectors from the corresponding bases (see Proposition \ref{DCFact_prop} below). This fact is crucial 
for an easy transfer of the proof of linear independence from the level $1$ case to higher levels.
In this sense, this proof is different than the proof given in \cite{P4}. In the $D_4$-case (cf. \cite{Ba}) it seems that this property does not hold
and it remains to see what would be a good way to capture phenomenons that are happening there. 
 
We give a brief outline of the paper. In sections \ref{S: Afine} and \ref{S: FS} we introduce the setting and pose the main 
problem.
In section \ref{S: ICDC} we find relations between monomials and use them to reduce the spanning set in terms of difference and initial conditions.  
The existence and some properties of intertwining operators and a simple current operator are established in sections \ref{S: IntOp}
and \ref{S: SimpCurr}. 
In section \ref{S: Rel} we explore the action of monomials of higher degree on the vectors from the top of Feigin-Stoyanovsky's type
subspaces. This will enable us to find suitable coefficients of intertwining operators having properties that we have discussed above,
and prove linear independence in the final section.   

We are grateful to Dra\v zen Adamovi\' c and Alex Feingold for useful informations and stimulating discussions on fusion rules and intertwining operators.

\section{Affine Lie algebra  $C_\ell^{(1)}$} \label{S: Afine}

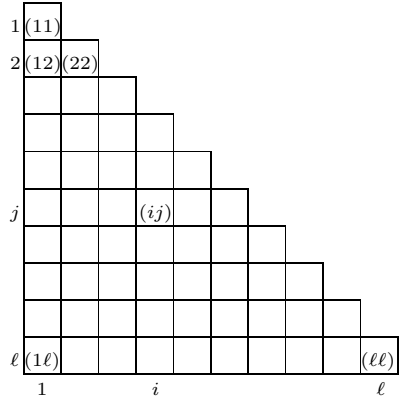
\begin{figure}[ht] \caption{The set of colors $\Gamma$} \label{Gamma_fig}
	\begin{center}\begin{picture}(150,150)(-10,-10) 
		\put(0,0){\line(1,0){140}} 
		\put(0,0){\line(0,1){140}}
		\put(140,0){\line(0,1){14}} 
		\put(0,140){\line(1,0){14}}
		\put(-5,129){$\scriptstyle 1$}
		\put(-5,115){$\scriptstyle 2$}
		\put(-5,3){$\scriptstyle \ell$} 
		\put(5,-8){$\scriptstyle 1$}
		\put(132,-8){$\scriptstyle \ell$}
		\linethickness{.075mm} 
		\put(0,14){\line(1,0){140}}
		\put(0,28){\line(1,0){126}}
		\put(0,42){\line(1,0){112}}
		\put(0,56){\line(1,0){98}}
		\put(0,70){\line(1,0){84}}
		\put(0,84){\line(1,0){70}}
		\put(0,98){\line(1,0){56}}
		\put(0,112){\line(1,0){42}}
		\put(0,126){\line(1,0){28}}
		\put(14,0){\line(0,1){140}}
		\put(28,0){\line(0,1){126}}
		\put(42,0){\line(0,1){112}}
		\put(56,0){\line(0,1){98}}
		\put(70,0){\line(0,1){84}}
		\put(84,0){\line(0,1){70}}
		\put(98,0){\line(0,1){56}}
		\put(112,0){\line(0,1){42}}
		\put(126,0){\line(0,1){28}}
		\put(-5,59){$\scriptstyle j$}
		\put(48,-8){$\scriptstyle i$}
		\put(0,129){$\scriptstyle (11)$}
		\put(0,115){$\scriptstyle (12)$}
		\put(14,115){$\scriptstyle (22)$}
		\put(0,3){$\scriptstyle (1 \ell)$}
		\put(126,3){$\scriptstyle (\ell \ell)$}
		\put(43,59){$\scriptstyle (ij)$}
		\end{picture}\end{center}
\end{figure}

Let ${\mathfrak g}$ be a complex simple  Lie algebra of type $C_\ell$ and let $\h$ be a Cartan
subalgebra of ${\mathfrak g}$. Let ${\mathfrak g}={\mathfrak h}+\sum {\mathfrak g}_\alpha$ be a root space decomposition of ${\mathfrak g}$. The corresponding root system $R$ may be realized in $\mathbb R\sp\ell$ with
the canonical basis $\epsilon_1,\dots,\epsilon_\ell$ as
$$
R=\{\pm \epsilon_i\pm \epsilon_j\,|\, 1\leq i\leq j\leq \ell\}\backslash\{0\}.
$$ 
We fix simple roots
$$\alpha_1=\epsilon_1-\epsilon_2,\quad\dots,\quad\alpha_{\ell-1}=\epsilon_{\ell-1}-\epsilon_\ell,\quad\alpha_\ell=2\epsilon_\ell$$ 
 and let   $ \g=\n_-+ \h + \n_+$ be the corresponding triangular decomposition. Let
$\theta=2\alpha_1+\dots+2\alpha_{\ell-1}+\alpha_\ell=2\epsilon_1$ be the maximal root and
$$
\omega_r= \epsilon_1+\dots+ \epsilon_r, \qquad r=1,\dots,\ell
$$ 
the fundamental weights. We fix root vectors $x_\alpha\in\g_\alpha$ and denote by $\alpha\sp\vee\in\mathfrak h$
dual roots. We identify $\mathfrak h$ and $\mathfrak h\sp*$ via the Killing form $\langle\,,\,\rangle$ normalized in such a way that  $\langle\theta,\theta\rangle=2$. We fix 
$$
\omega=\omega_\ell=\epsilon_1+\dots+\epsilon_\ell \in\h^*.
$$ 
This is the minuscule coweight, that is 
$$\langle\omega,\alpha\rangle\in\{-1,0,1\}\quad\text{for all}\quad \alpha\in R,
$$
and hence we have a $\mathbb Z$-gradation of $\g$
$$
 \mathfrak g =\mathfrak g_{-1} + \mathfrak g_0 + \mathfrak g_1 , \qquad
 {\mathfrak g}_0 = {\mathfrak h} +
\sum_{\langle\omega,\alpha\rangle=0}\, {\mathfrak g}_\alpha,
\qquad
\displaystyle {\mathfrak g}_{\pm1}  = 
\sum_{\langle\omega,\alpha\rangle= \pm 1}\, {\mathfrak g}_\alpha.
$$
Note that
\begin{equation}
\Gamma  =
\{\,\alpha \in R \mid \langle\omega,\alpha\rangle = 1\}  =  \{\epsilon_i+\epsilon_j\,|\, 1\leq i \leq j\leq \ell\}.
\end{equation}
We say that $\Gamma$ is {\it the set of colors} and we write
\begin{equation}
(ij)=\epsilon_i+\epsilon_j\in\Gamma\quad\text{and}\quad x_{ij}=x_{\epsilon_i+\epsilon_j}
\end{equation}
(see figure \ref{Gamma_fig}). 

The subspaces ${\mathfrak g}_{\pm1}\subset\mathfrak g$ are commutative subalgebras, $\g_0$ is reductive and $[\g_0,\g_0]$
is a simple algebra of type $A_{\ell-1}$ with root basis $\alpha_1,\dots,\alpha_{\ell-1}$.
We identify $[\g_0,\g_0]$ with the Lie algebra $\mathfrak{sl}(\ell,\mathbb C)$ acting on 
the canonical basis $e_1,\dots,e_\ell$ of the vector space $\C^\ell$ by the rule 
\begin{equation} \label{VectRep_eq}
e_i \xrightarrow{x_{-\alpha_i}} e_{i+1}, 
\end{equation}
The subalgebra $\g_0$ acts on $\g_1$ by adjoint action. For a suitably chosen root vectors $x_{ij}$ 
this action is given by 
 \begin{eqnarray}
\label{E:g_0 action1} x_{ij} & \xrightarrow{x_{-\alpha_i}} & x_{i+1,j}\quad \textrm{for}\ i<j, \\
 x_{ij} & \xrightarrow{x_{-\alpha_j}} & x_{i,j+1}\quad \textrm{for}\ i<j, \\
\label{E:g_0 action2} x_{ii} & \xrightarrow{x_{-\alpha_i}} & 2x_{i,i+1}
 \end{eqnarray}
(cf. \cite{H}).

We identify the Weyl group $\Weyl$ of $[\g_0,\g_0]$ with the group of permutations
$$
\sigma\colon i \mapsto \sigma(i), \quad i=1,\dots,\ell,
$$
so that for $\alpha=\epsilon_{i}-\epsilon_{j}$, $i\neq j$, we have 
$\sigma\alpha=\epsilon_{\sigma(i)}-\epsilon_{\sigma(j)}$. For each $\sigma\in\Weyl$ we have an
automorphism $\sigma$ of $\g_0$ and a linear map $\sigma$ on the vector representation $V_1=\C^\ell$
of $\g_0$,
\begin{equation} \label{E: Weyl map}
\sigma\colon \g_0\to \g_0,\qquad \sigma\colon V_1\to V_1,
\end{equation}
such that
\begin{equation} \label{E: TildaNajava1}
\sigma\,x_\alpha\in\mathbb C\sp\times x_{\sigma\alpha},\qquad \sigma\,e_i\in\mathbb C\sp\times e_{\sigma(i)}
\end{equation}
and
\begin{equation}\label{E: action of sigma 1}
\sigma (xv)=(\sigma x)(\sigma  v)
\end{equation}
for $x\in \g_0$ and $v\in V_1$. For simple reflection $\sigma_i\in\Weyl$ the linear map 
$\sigma_i$ is $(\exp x_{-\alpha_i})(\exp -x_{\alpha_i})(\exp x_{-\alpha_i})$ 
and formula (\ref{E: action of sigma 1}) holds in general for integrable
$\g_0$-modules (cf. \cite{K}).
Since $\g_1\cong S^2(\C^\ell)$, we also have a linear map $\sigma$ on $\g_1$ such that
\begin{equation} \label{E: TildaNajava2}
 \sigma\,x_{ij}\in\mathbb C\sp\times x_{\sigma(i)\sigma(j)}.
\end{equation}

To abbreviate expressions like the ones in \eqref{E: TildaNajava1} and \eqref{E: TildaNajava2} we introduce the following
notation: for two vectors, monomials, etc. we write $x \sim y$ if
the two are equal up to a nonzero scalar, i.e.
\begin{equation}\label{E: TildaDef}
x = Cy, \qquad \textrm{for some}\ C\in\C^\times.
\end{equation}
In this way,  relations \eqref{E: TildaNajava1} and \eqref{E: TildaNajava2} can be rewritten in the following way:
\begin{eqnarray*} 
 & \sigma\,x_\alpha \sim  x_{\sigma\alpha},\qquad \sigma\,e_i \sim e_{\sigma(i)}, & \\
 & \sigma\,x_{ij} \sim x_{\sigma(i)\sigma(j)}. &
\end{eqnarray*}

\bigskip


\smallskip
Denote by $\tilde{\mathfrak g}$ the affine Lie algebra of type $C_\ell\sp{(1)}$ associated  to $\g$,
 $$
 \hat{\mathfrak g} = \mathfrak{g} \otimes
 \mathbb{C}[t,t^{-1}] + \mathbb{C}c, \qquad
  \tilde{\mathfrak g} = \hat{\mathfrak g} + \mathbb{C} d,
 $$
with the canonical central element $c$ and the degree element $d$ . Set
$$
x(n)=x\otimes t^{n}
$$  
for $x\in{\mathfrak g}$ and $n\in\mathbb Z$ and denote by  $ x(z)=\sum_{n\in\mathbb Z}
 x(n) z^{-n-1}$ a formal Laurent series  in formal variable $z$. The commutation relations in 
 $ \tilde{\mathfrak g}$ are
 $$
  [x(i),y(j)]= [x,y](i+j)+ i\langle x,y\rangle
 \delta_{i+j,0}c,\qquad [c,\gt]  =  0, \qquad [d,x(j)]  =  j x(j). 
 $$
We have a triangular decomposition 
$$
\gt=\tilde{\n}_-+ \tilde{\h} + \tilde{\n}_+,
$$
where
$$
\tilde{\n}_-=\n_-+ \g\otimes t^{-1}\C [t^{-1}],\qquad
\tilde{\h}=\h  + \mathbb{C}c + \mathbb{C} d,\qquad
\tilde{\n}_+=\n_++ \g\otimes t\C [t].
$$
We also have the induced $\Z$-gradation 
$$
\gt=\gt_{- 1}+\gt_{0}+\gt_{ 1}
$$
of  affine Lie algebra $\gt$, where
$$
\gt_0 = {\mathfrak g}_0\otimes\C [t,t^{-1}]\oplus \C c \oplus \C d,\qquad
\gt_{\pm 1} = {\mathfrak g}_{\pm 1}\otimes\C [t,t^{-1}].
$$
The subspace $\gt_1\subset \gt$ is a commutative subalgebra and $\g_0$ acts on
$\gt_1$ by adjoint action.

We denote by $\Lambda_0,\dots,\Lambda_\ell$ the fundamental weights of $\gt$,
\begin{equation} \label{TopSl2_rel}
\Lambda_r=\Lambda_0+\omega_r,\quad r=1,\dots,\ell.
\end{equation}


\section{Feigin-Stoyanovsky's type subspaces} \label{S: FS}

Denote by $L(\Lambda)$ a standard (i.e. integrable highest weight) $\tilde{\mathfrak g}$-module with the highest weight
$$
\Lambda=k_0 \Lambda_0+k_1 \Lambda_1+\dots+k_\ell \Lambda_\ell,
$$
$k_i\in\Z_+$ for $i=0,\dots,\ell$. Throughout the paper we denote by
$k=\Lambda(c)$ the level of  $\tilde{\mathfrak g}$-module $L(\Lambda)$,
$$
k=k_0 +k_1 +\dots+k_\ell ,
$$
and by $v_\Lambda$ a fixed highest weight vector of $L(\Lambda)$.

For each integral dominant $\Lambda$ we have a  {\it Feigin-Stoyanovsky's type subspace}
$$
W(\Lambda)=U(\widetilde{{\mathfrak
g}}_1)v_\Lambda\subset L(\Lambda).
$$
This space has a Poincare-Birkhoff-Witt spanning set 
\begin{equation} \label{E: PBW base def}
\{x_{\gamma_1}(-n_1)
x_{\gamma_2}(-n_2)\cdots x_{\gamma_t}(-n_t) v_\Lambda\,|\,
t\in\Z_+,\gamma_i\in \Gamma, n_i\in\N\}.
\end{equation}
The main problem in this paper is to reduce this PBW spanning set to a basis of
$W(\Lambda)$. The main three steps in our construction are:
\begin{itemize}
 \item find relations for standard modules,
\item reduce the spanning set, and
\item prove linear independence by
using intertwining operators.
 \end{itemize}

\section{Difference and initial conditions} \label{S: ICDC}

Start from the vertex operator algebra relation (cf. \cite{LP}, \cite{MP}, \cite{LL})
\begin{equation} \label{FKS_rel}
x_\theta(z)^{k+1}=\sum_{N\in\Z}\left( \sum_{n_1+\dots+n_{k+1}=N} x_{\theta}(-n_1)\cdots x_\theta (-n_{k+1})\right)z^{n-k-1}=0\qquad \textrm{on} \ L(\Lambda).
\end{equation}
Adjoint $\g_0$-action on \eqref{FKS_rel} gives us the space of relations $U(\g_0)\cdot x_\theta(z)^{k+1}=0$.
This is a finite-dimensional $\g_0$-module with the highest weight $2(k+1)\omega_1$.
Hence, as a vector space, it is isomorphic to $S^{2(k+1)}(\C^\ell)$. The basis of $S^{2(k+1)}(\C^\ell)$ is given by
$e_1^{m_1}\cdots e_\ell^{m_\ell}$, $m_1+\dots+m_\ell=2(k+1)$, which we view as multisets
$\{1^{m_1},\dots,\ell^{m_\ell}\}$. Since relations 
\eqref{E:g_0 action1}--\eqref{E:g_0 action2} hold, one can easily see
that the corresponding ``basis'' of the set of relations is given by the following proposition:

\begin{prop}
On $L(\Lambda)$ the following relations hold	
\begin{equation}\label{DC_Frel}
\sum_{\{i_1,\dots,i_{k+1}\}\cup\{ j_1,\dots,j_{k+1} \}=\{1^{m_1},\dots,\ell^{m_\ell}\}  } \hspace{-3ex} C_{\bf ij} x_{i_1 j_1}(z)x_{i_2 j_2}(z)
\cdots x_{i_{k+1} j_{k+1}}(z)=0,
\end{equation} 
for some $C_{\bf ij}\in\C^\times$, where the sum runs over all such partitions of the multiset $\{1^{m_1},\dots,\ell^{m_\ell}\}$. 
\end{prop}

For each power of $z$ from \eqref{DC_Frel} we obtain a relation between monomials
\begin{equation}
\label{DC_rel}
\sum_{\substack{n_1+\dots+n_{k+1}=N\\ \{i_1,\dots,i_{k+1},j_1,\dots,j_{k+1} \}=\{1^{m_1},\dots,\ell^{n_\ell}\} }} 
\hspace{-3ex} C_{\bf ij} x_{i_1 j_1}(-n_1)x_{i_2 j_2}(-n_2)\cdots x_{i_{k+1} j_{k+1}}(-n_{k+1})=0.
\end{equation}
We find the smallest monomials in these relation,
the {\em leading terms} of relations. Since they can be expressed 
as a sum of higher terms, we can exclude them from the spanning set \eqref{E: PBW base def}.

We introduce a linear order on monomials in the following way.
First, define a linear order on the set of colors $\Gamma$:
$(i'j')<(ij)$ if $i'>i$ or $ i'=i,\, j'>j$. 
On the {\em set of variables} $\Gamt=\{x_\gamma(n) \,|\,\gamma\in\Gamma, n\in \Z \}$  define a linear order by
$x_\alpha(n)<x_\beta(n')$ if $n<n'$ or
$n=n',\, \alpha<\beta$. 
Assume that variables in monomials are sorted descendinly from right
to left. The order $<$ on the set of monomials is
defined as a lexicographic order, where we compare variables from
right to left (from the greatest to the lowest one).

Order $<$ is compatible with multiplication (see [P1], [T1]):
\begin{equation} \label{OrdMult_rel}
\textrm{if} \quad x(\pi)<x(\pi')\quad \textrm{then}\quad x(\pi)x(\pi_1)<x(\pi')x(\pi_1)\quad 
\end{equation}
for monomials $x(\pi),x(\pi'),x(\pi_1)\in\C[\Gamt]$.

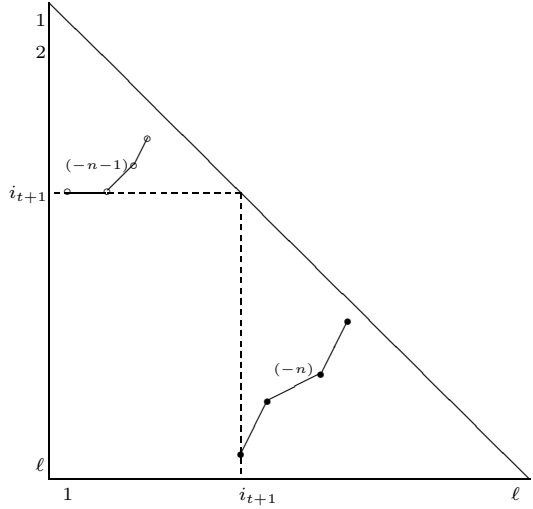
\begin{figure}[ht] \caption{Difference conditions for level $k$} \label{DCSeq_fig}
\begin{center}\begin{picture}(200,200)(-10,-10) 
\put(0,0){\line(1,0){180}} \put(0,0){\line(0,1){180}}
\put(180,0){\line(-1,1){180}} 
\put(-5,171){$\scriptstyle 1$}
\put(-5,159){$\scriptstyle 2$}
\put(-5,3){$\scriptstyle \ell$} \put(5,-8){$\scriptstyle
1$}
\put(173,-8){$\scriptstyle \ell$}
\put(84,40){$\scriptscriptstyle (-n)$}
\put(110,58){$\scriptscriptstyle \bullet$}
 \put(102,40){\line(1,2){10}}
\put(100,38){$\scriptscriptstyle \bullet$}
 \put(82,30){\line(2,1){20}}
\put(80,28){$\scriptscriptstyle \bullet$}
 \put(72,10){\line(1,2){10}} 
\put(70,8){$\scriptscriptstyle  \bullet$}
 \multiput(72,108)(-4,0){18}{\line(-1,0){2}}
\multiput(72,108)(0,-4){27}{\line(0,-1){2}}
\put(-15,106){$\scriptstyle i_{t+1}$}
\put(71,-8){$\scriptstyle i_{t+1}$}

\put(35,127){$\scriptscriptstyle \circ$}
\put(30,117){$\scriptscriptstyle \circ$}
\put(20,107){$\scriptscriptstyle \circ$}
\put(5,107){$\scriptscriptstyle \circ$}
\put(7,108){\line(1,0){14}} 
\put(32,119){\line(1,2){5}} 
\put(22,109){\line(1,1){10}}
\put(6,117){$\scriptscriptstyle
	(-n-1)$}

\end{picture}\end{center}
\end{figure}

The leading terms can be most conveniently described in terms of exponents: 

\begin{prop}
A monomial
$$ x(\pi') =  x_{i_t
j_t}(-n-1)^{b_{i_t j_t}}\cdots x_{i_1 j_1}(-n-1)^{b_{i_1 j_1}}  x_{i_{s} j_{s}}(-n)^{a_{i_{s} j_{s}}}\cdots x_{i_{t+1}
j_{t+1}}(-n)^{a_{i_{t+1} j_{t+1}}},$$
where $i_1 \leq \dots \leq i_t
\leq j_t\leq\ldots\leq j_1
\leq i_{t+1} \leq \dots \leq i_{s} \leq j_{s} \leq \ldots\leq j_{t+1}$, $(i_\nu,j_\nu)\neq (i_{\nu+1},j_{\nu+1})$, and
\begin{equation}
\label{LT_eq}
b_{i_1 j_1}+\ldots+b_{i_t
j_t}+a_{i_{t+1} j_{t+1}}+\ldots+a_{i_s j_s}=k+1,
\end{equation}
 is a leading term of a relation \eqref{DC_rel} corresponding to the multiset 
 \begin{eqnarray*}
 & & \hspace{-3ex} \{{i_1}^{b_{i_1 j_1}},\dots,{i_t}^{b_{i_t j_t}},\ 
 {j_t}^{b_{i_t j_t}} ,\dots, {j_1}^{b_{i_1 j_1}},\\
 & & \hspace{10ex}  
 {i_{t+1}}^{a_{i_{t+1} j_{t+1}}},\dots, {i_{s}}^{a_{i_{s} j_{s}}}, 
 {j_{s}}^{a_{i_{s} j_{s}}},\dots,{j_{t+1}}^{a_{i_{t+1} j_{t+1}}}\}
 \end{eqnarray*}
and a degree $N=(k+1)n + m$, where $m=b_{i_1 j_1}+\ldots+b_{i_t j_t}$.
\end{prop}
The colors of leading terms lie on diagonal paths in $\Gamma$, see figure \ref{DCSeq_fig}.

\begin{dokaz}
	Consider a relation \eqref{DC_rel} corresponding to a multiset $\{p_1\leq\dots\leq p_{2k+2}\}$ and
	a total degree $-n_1-\dots-n_{k+1}=-N\in -\N$.
	
	First consider the case $N=(k+1)n$. In this case the leading term has all factors of the same degree $-n$
	so we need to find minimal configuration of colors whose rows and columns joined give $\{p_1,\dots,p_{2k+2}\}$.
	It is clear that rows of the minimal configuration are $\{p_1,\dots,p_{k+1}\}$ and its columns are 
	$\{p_{k+2},\dots,p_{2k+2}\}$. Otherwise there would be a leading term
	$x(\pi)=x(\pi_1)x_{ij}(-n)x_{i'j'}(-n)$, $i\leq j< i'\leq j'$, and it is clear that
	a monomial $x(\pi')=x(\pi_1)x_{ij'}(-n)x_{ji'}(-n)$ from the same relation 
	is smaller than $x(\pi)$. By a similar argument we conclude that the minimal configuration is
	obtained by pairing maximal rows with minimal columns; i.e. it consists of colors $\{(p_1 p_{2k+2}),
	(p_2 p_{2k+1}), \dots, (p_{k+1} p_{k+2})\}$ (if $x(\pi)=x(\pi_1)x_{ij}(-n)x_{i'j'}(-n)$, $i< i'< j < j'$, then 
	$x(\pi')=x(\pi_1)x_{ij'}(-n)x_{i'j}(-n)$ is smaller than $x(\pi)$). Hence it is a configuration whose colors lie on a diagonal path
	like in figure \ref{DCSeq_fig}.
	
	Next, consider the case $N=(k+1)n+m$. In this case the leading term has $m$ factors of degree $-n-1$, and 
	$k+1-m$ factors of degree $-n$. The leading term is obtained first by choosing the minimal possible $(-n)$-part, 
	and then the minimal possible $(-n-1)$-part. Hence $(-n)$-part corresponds to $\{p_{2m+1,\dots,p_{2k+2}}\}$, while
	$(-n-1)$-part corresponds to $\{p_1,\dots, p_{2m}\}$, and colors of these parts lie on diagonal paths 
	like in figure \ref{DCSeq_fig}.
\end{dokaz}

We say that a monomial  $x(\pi)$ satisfies {\em
difference conditions}, or shortly, that $x(\pi)$
{\em satisfies $DC$}, if it does not contain leading terms. 
More precisely, $x(\pi)$ satisfies difference conditions if for
any $n\in\N$ and $i_1 \leq \dots \leq i_t
\leq j_t\leq\ldots\leq j_1
\leq i_{t+1} \leq \dots \leq i_{s} \leq j_{s} \leq \ldots\leq j_{t+1}$,
\begin{equation} \label{DC_def}
b_{i_1 j_1}+\ldots+b_{i_t j_t}+a_{i_{t+1} j_{t+1}}+\ldots+a_{i_s
j_s}\leq k,
\end{equation}
where $a_{ij}$'s and $b_{ij}$'s denote exponents of $x_{ij}(-n)$ and
$x_{ij}(-n-1)$ in $x(\pi)$, respectively.

Note that in the case of level $k=1$, difference conditions imply
that if  $x(\pi)=x(\pi') x_{ij}(-n)$ then $x(\pi')$ does not contain
factors $x_{i' j'}(-n)$, $i\leq i'\leq  j'\leq j$ or $i'\leq i\leq j\leq j'$,
nor it contains factors   
$x_{i' j'}(-n-1)$, $i'\leq j' \leq i$. Hence, 
$x(\pi)=\dots x_{i_s'
j_s'}(-n-1)\cdots x_{i_1' j_1'}(-n-1)  x_{i_{t} j_{t}}(-n) \cdots x_{i_1
j_1}(-n)\dots$
satisfies difference conditions for level $k=1$ if
\begin{equation} \label{DClev1_eq}
i_1<\dots <i_t,\quad j_1<\dots<j_t,\quad  i_1'<\dots < i_s',\quad i_t<j_1'<\dots<j_s'. 
\end{equation}
Its colors lie on diagonal paths and a diagonal path of 
$(-n-1)$-part lies below $i_t$-th row, where $i_t$ is the column of the smallest
color of the $(-n)$-part; see figure \ref{DClev1_fig}.

\begin{napomena}
	{\em
		Similar difference conditions appear in another construction of combinatorial bases for $C_\ell^{(1)}$ (cf. \cite{PS}).
	}\end{napomena}

\begin{lm} \label{IC1_lm}
On $L(\Lambda_r)$
\begin{eqnarray*}
x_{ij}(-1)v_{\Lambda_r}=0, & & \textrm{if} \  j\leq r,\\
x_{ij}(-1)v_{\Lambda_r}\neq 0, & & \textrm{if} \  j> r.
\end{eqnarray*}

\end{lm}

\begin{dokaz}
For $\alpha\in R$ denote by $\gsl_2 (\alpha)\subset \g$ a subalgebra generated by $x_{\alpha}$ and
$x_{-\alpha}$, and let
$$\tilde{\gsl}_2(\alpha)=\gsl_2 (\alpha)\otimes \C[t,t^{-1}]\oplus \C c \oplus \C d\subset \gt$$
be the corresponding affine Lie algebra of type $A_1^{(1)}$. It has a canonical central
element $c'=\langle x_{\alpha}, x_{-\alpha}\rangle c=2c/\langle \alpha,\alpha \rangle$. Hence
the restriction of $L(\Lambda_r)$ is a level one representation if $\alpha$ is a long root, and
a level two representation if $\alpha$ is a short root. 

Consider $\alpha=(jj)=2\epsilon_j$, $j\leq r$. Then $(jj)^\vee=(jj)=2\epsilon_j$ and $\langle \omega_r,(jj)^\vee\rangle = \delta_{j>r}$. 
Hence, by \eqref{TopSl2_rel}, if $j\leq r$, then
$U(\tilde{\gsl}_2(\alpha))v_{\Lambda_r}$ is a level one representation with
$2$-dimensional $\gsl_2(\alpha)$-module on top, and therefore it must be
the standard $A_1^{(1)}$-module $L(\Lambda_1)$. If $j> r$, then the $\gsl_2(\alpha)$-module on top is $1$-dimensional, 
hence $U(\tilde{\gsl}_2(\alpha))v_{\Lambda_r}$ must be the standard $A_1^{(1)}$-module $L(\Lambda_0)$ (cf. \cite{K}). Therefore
\begin{eqnarray}
\label{ICLong_eq1} x_{jj}(-1)v_{\Lambda_r}=0, & & \textrm{if} \  j\leq r,\\
\label{ICLong_eq2} x_{jj}(-1)v_{\Lambda_r}\neq 0, & & \textrm{if} \  j> r,
\end{eqnarray}
which proves the lemma for $i=j$.

For $i<j$, action by $x_{\alpha_i}(0)\cdots x_{\alpha_{j-1}}(0)$ 
on \eqref{ICLong_eq1} and \eqref{ICLong_eq2} gives the claim. 
\end{dokaz}

\begin{figure}[ht] \caption{Initial conditions for $W(\Lambda_r)$} \label{IClev1_fig}
\begin{center}\begin{picture}(200,200)(-10,-10) 
\put(0,0){\line(1,0){180}} \put(0,0){\line(0,1){180}}
\put(180,0){\line(-1,1){180}} 
\put(-5,171){$\scriptstyle 1$}
\put(-5,159){$\scriptstyle 2$}
\put(-5,3){$\scriptstyle \ell$} \put(5,-8){$\scriptstyle
1$}
\put(173,-8){$\scriptstyle \ell$}
\put(83,50){$\scriptscriptstyle (-1)$}
\put(60,70){$\scriptscriptstyle \bullet$}
\put(70,50){$\scriptscriptstyle \bullet$}
\put(90,40){$\scriptscriptstyle \bullet$}
\put(120,30){$\scriptscriptstyle  \bullet$}
\put(62,72){\line(1,-2){10}}
\put(72,52){\line(2,-1){20}}
\put(92,42){\line(3,-1){30}}
\multiput(100,80)(-4,0){25}{\line(-1,0){2}}
\put(-5,78){$\scriptstyle r$}
\end{picture}\end{center}
\end{figure}
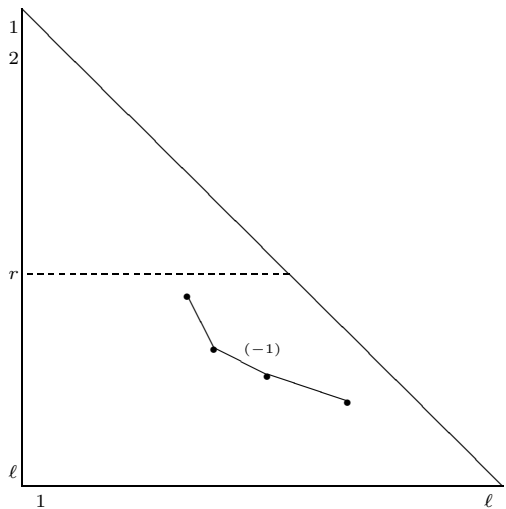

A monomial $x(\pi)$ satisfies {\em initial conditions} for $W(\Lambda_r)$ if 
it does not contain a factor $x_{ij}(-1)$, $i\leq j \leq r$. 
Note that if a monomial $x(\pi)$ satisfies difference and initial conditions
for $W(\Lambda_r)$, then the colors of $(-1)$-factors lie on a diagonal path  
bellow the $r$-th row (see figure \ref{IClev1_fig}).

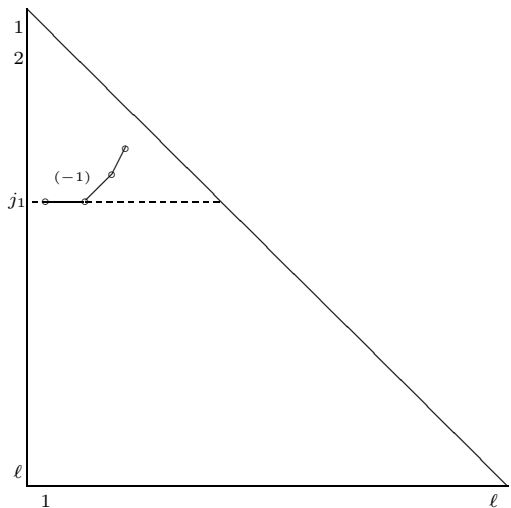
\begin{figure}[ht] \caption{Initial conditions for $\Lambda=k_0\Lambda_0+\ldots+k_\ell\Lambda_\ell$} \label{ICseq_fig}
\begin{center}\begin{picture}(200,200)(-10,-10) 
\put(0,0){\line(1,0){180}} \put(0,0){\line(0,1){180}}
\put(180,0){\line(-1,1){180}} 
\put(-5,171){$\scriptstyle 1$}
\put(-5,159){$\scriptstyle 2$}
\put(-5,3){$\scriptstyle \ell$} \put(5,-8){$\scriptstyle
1$}
\put(173,-8){$\scriptstyle \ell$}
\put(35,126){$\scriptscriptstyle \circ$}
\put(30,116){$\scriptscriptstyle \circ$}
\put(20,106){$\scriptscriptstyle \circ$}
\put(5,106){$\scriptscriptstyle \circ$}
\put(8,107){\line(1,0){13}} 
\put(32,118){\line(1,2){5}} 
\put(22,108){\line(1,1){10}}
\put(10,115){$\scriptscriptstyle
(-1)$}
\multiput(72,107)(-4,0){18}{\line(-1,0){2}}
\put(-7,106){$\scriptstyle j_1$}
\end{picture}\end{center}
\end{figure}

Generally, let $\Lambda=k_0\Lambda_0+\ldots+k_\ell\Lambda_\ell$. We say that
$x(\pi)$ satisfies {\em initial conditions} for $W(\Lambda)$ if 
for every $i_1 \leq \dots \leq i_t
\leq j_t\leq\ldots\leq j_1$,
\begin{equation} \label{IC_def}
a_{i_1 j_1}+\ldots+a_{i_t
j_t}\leq k_0+k_1+\ldots + k_{j_1-1}
\end{equation}
where $a_{ij}$'s denote exponents of $x_{ij}(-1)$ in $x(\pi)$
(see figure \ref{ICseq_fig}).
One immediately sees that for $\Lambda=\Lambda_r$ the two definitions of initial
conditions are equivalent.

\begin{napomena}\label{ICred_nap}
{\em
Like in \cite{Ba} and \cite{T2}, initial conditions can be expressed in terms of difference
conditions by adding ``imaginary'' $(0)$-factors to $x(\pi)$. Let
$$x(\pi_0)=x_{1\ell}(0)^{k_1} x_{2\ell}(0)^{k_2}\cdots x_{\ell \ell}(0)^{k_\ell}.$$ 
Note that colors of $x(\pi_0)$ lie on a diagonal path like in figure \ref{DCSeq_fig}.
Then $x(\pi)$ satisfies difference and initial conditions for $W(\Lambda)$ if
and only if $x(\pi')=x(\pi)x(\pi_0)$ satisfies difference conditions. In
fact, initial conditions are defined in such a way so that this property holds.
}\end{napomena}

\begin{prop}
The set of monomial vectors $x(\pi) v_\Lambda$ satisfying difference conditions
\eqref{DC_def} and initial conditions \eqref{IC_def} span $W(\Lambda)$. 
\end{prop}

\begin{dokaz}
By \eqref{DC_rel}, if $x(\pi)$ does not satisfy difference condition
\eqref{DC_def}, then $x(\pi)v_\Lambda$ can be expressed in terms of higher
monomial vectors. Hence we can exclude $x(\pi)v_\Lambda$ from the spanning set \eqref{E: PBW base def}.

By definition of level one initial conditions, if $x(\pi)$ does not satisfy initial conditions
for $W(\Lambda_r)$, then $x(\pi)v_{\Lambda_r}=0$ and $x(\pi)v_{\Lambda_r}$ can be excluded from
the spanning set \eqref{E: PBW base def}.   
Assume $k>1$ and that $x(\pi)$ does not satisfy initial condition 
\eqref{IC_def}, and set
$d=k_0+k_1+\ldots + k_{j_1-1}+1$.  Factorize $x(\pi)=x(\pi'')x(\pi')$, where
$x(\pi')$ consists only of  $(-1)$-factors lying on a diagonal path from \eqref{IC_def}. 
Furthermore, one
can assume that the length of $x(\pi')$ is equal to $d$ and that $d<k+1$ (otherwise \eqref{IC_def}
is equivalent to \eqref{DC_def}). 



Set
$\Lambda' = \sum_{r=0}^{j_1-1} k_r \Lambda_r$, $\Lambda'' = \sum_{r=j_1}^{\ell} k_r \Lambda_r=\Lambda - \Lambda'$.
Denote by $v_{\Lambda'}$ and $v_{\Lambda''}$ highest weight vectors
of  standard modules $L(\Lambda')$ and $L(\Lambda'')$. Then, by
the complete reducibility, $L(\Lambda)  \subset  L(\Lambda')\otimes L(\Lambda'')$,
$v_\Lambda = v_{\Lambda'}\otimes v_{\Lambda''}$.
Since, by lemma \ref{IC1_lm}, all factors of $x(\pi')$ annihilate $v_{\Lambda''}$, we have
$$x(\pi')v_\Lambda=(x(\pi')v_{\Lambda'})\otimes v_{\Lambda''}.$$
Note that $L(\Lambda')$ is a module of level $d-1<k$.
From relations \eqref{DC_rel} for the module $L(\Lambda')$ we obtain
monomials $x(\pi_1'),\dots,x(\pi_s')$ 
such that $x(\pi')v_{\Lambda'}=C_1
x(\pi_1')v_{\Lambda'}+\dots+C_s x(\pi_s')v_{\Lambda'}$,
$C_t\in\C^\times$, and $x(\pi')<x(\pi_t')$. Also from these relations,
we see that colors of monomials $x(\pi_t')$ lie in the $j_1$-th row or above. 
By lemma \ref{IC1_lm}, all factors of $x(\pi_t')$ also act as $0$ on $v_{\Lambda''}$. Consequently,
$$x(\pi')v_{\Lambda}=C_1 x(\pi_1')v_{\Lambda}+\dots+C_s
x(\pi_s')v_{\Lambda},$$
and $x(\pi)v_\Lambda$ can be expressed in terms of higher
monomial vectors. Therefore, it can be excluded from the spanning set \eqref{E: PBW base def}.
\end{dokaz}

Like in the  $A_\ell$-case ([T2]), difference and initial conditions for level $k>1$ can be interpreted in
terms of difference and initial conditions for level $1$:

\begin{prop}
\label{DCFact_prop} Let $L(\Lambda)\subset L(\Lambda_{i_1})\otimes \cdots \otimes L(\Lambda_{i_k})$ be a standard module of level $k$. Monomial
$x(\pi)$ satisfies difference and initial conditions
for $W(\Lambda)$ if and only if there exists a factorization
 $$x(\pi)=x(\pi^{(1)})\cdots x(\pi^{(k)}),$$
such that  $x(\pi^{(j)}) $ satisfies difference and initial
conditions for $W(\Lambda_{i_j})$.
\end{prop}

\begin{dokaz}
We follow the idea from [T2]; here we give a sketch of the proof. 
First consider the case $\Lambda=k\Lambda_0$ for which initial conditions
do not provide any additional relations and  we only need to consider difference
conditions.

Define another order on the set of variables:
$x_{ij}(-n) \sqsubset x_{i'j'}(-n')$ if either $-n \leq -n'-2$
or $-n=-n'-1,j>i'$ or $-n=-n',i>i', j>j'$. 
Equivalently,  $x_{ij}(-n) \sqsubset x_{i'j'}(-n')$  if $x_{ij}(-n) < x_{i'j'}(-n')$ and
a monomial $x_{ij}(-n) x_{i'j'}(-n')$ satisfies level $1$ difference conditions. 
This is a strict partial order on $\Gamt$.

Consider monomials $x(\pi)\in \C[\Gamt]$ as multisets; then a monomial $x(\pi)$ satisfies 
level $k$ difference conditions  if and only if every subset of $x(\pi)$ in which there are
no two elements comparable in the sense of $\sqsubset$, has at most $k$ elements.
To see this, let $x_{ij}(-n),x_{i'j'}(-n')\in\Gamt$, $x_{ij}(-n)< x_{i'j'}(-n')$. 
They are incomparable in the sense $\sqsubset$ if and only if
either $-n=-n'-1$ and $j\leq i'$ or $-n=-n'$ and $i\geq i'$ or $j\geq j'$.
It is clear that factors of leading terms \eqref{LT_eq} are mutually
incomparable; consequently, if $x(\pi)$
does not satisfy difference conditions, then it has a subset of at
least $k+1$ mutually incomparable elements. Conversely, consider a
subset of $x(\pi)$ whose elements are mutually incomparable. 
By the observation above, degrees of its elements differ for at most $1$. 
Moreover, elements of the same degree lie on a diagonal path and the two 
paths are related like in \eqref{LT_eq}.

Note that if $x_{\gamma_1}(-n_1) \sqsubset \dots \sqsubset x_{\gamma_t}(-n_t)$, 
then the monomial $x_{\gamma_1}(-n_1)\cdots x_{\gamma_t}(-n_t)$ satisfies 
level $1$ difference conditions. Now, a combinatorial lemma from [T2]
implies that $x(\pi)$ can be partitioned into $k$ linearly ordered subsets
which proves proposition in the $\Lambda=k\Lambda_0$ case.

By remark \ref{ICred_nap}, in the general case $\Lambda=k_0\Lambda_0+\dots+k_\ell \Lambda_\ell$,
initial conditions can be regarded as difference conditions by considering monomials
$x(\pi')=x(\pi) x(\pi_0)$ instead of $x(\pi)$, where $x(\pi_0)=x_{1\ell}(0)^{k_1} \cdots x_{\ell \ell}(0)^{k_\ell}$. 
By above arguments, we can partition $x(\pi')$ into $k$ linearly ordered subsets. 
Since $x_{i\ell}(0)$'s are mutually incomparable, they lie in different subsets. Moreover,
a subset containing  $x_{i\ell}(0)$ gives a monomial satisfying difference and initial conditions for
$W(\Lambda_i)$, again by remark \ref{ICred_nap}.
\end{dokaz}

\section{Intertwining operators} \label{S: IntOp}

Consider a $\g_0$-module $V_i=U(\g_0)v_{\Lambda_i}\subset L(\Lambda_i)$ and $\C^\ell$ as the vector representation for $\g_0$-action. Then  
$$V_i=\bigwedge^i \C^\ell,$$
 for $i=0,\dots,\ell$ (cf. \cite{Bo}, \cite{H}).
If $e_1,\dots,e_\ell$ is a basis for $\C^\ell$ with $\g_0$-action
defined by \eqref{VectRep_eq}, then  
\begin{equation}
\label{Top_vect}
v_{p_1 \dots p_i}=e_{p_1}\wedge\cdots\wedge e_{p_i}, 
\qquad 1\leq p_1<\dots<p_i\leq\ell
\end{equation} 
is a basis for $V_i$.
Moreover, $$v_{\Lambda_i}=v_{1 2 \dots i}.$$ 
If $I=\{p_1,\dots,p_i\}$, $1\leq p_1<\dots<p_i\leq\ell$, then denote by $v_I=v_{p_1 \dots p_i}$. Note that for
each $\sigma\in\Weyl$ and each $i=1,\dots,\ell$ we have a linear map $\sigma$ on $V_i$ such that
$$
\sigma v_{p_1 \dots p_i} \sim v_{\sigma(p_1) \dots \sigma(p_i)},
$$
where we use the notation from \eqref{E: TildaDef}.
Later on we shall also use other consequences of the formula (\ref{E: action of sigma 1}) for integrable $\g_0$-modules, for example the formula
$$
\sigma (x_{pq}(-1)v_{p_1 \dots p_i})\sim  x_{\sigma(p)\sigma(q)}(-1)v_{\sigma(p_1) \dots \sigma(p_i)}.
$$
\bigskip

\begin{lm}\label{Top_rel2} Let $v\in V_i$. Then  
$x_\gamma(n)v=0$, 
for $\gamma\in\Gamma$, $n\geq 0$, and 
$x_\alpha(n)v=0$, 
for $\alpha\in R$, $n\geq 1$.\end{lm}
\begin{dokaz}
We only need to show $xv=0$ for all $x\in\g_1, v\in V_i$; other relations are clear from the definition.
Let $v'\in V_i$ be such that $\g_1 v'=0$. Let $x\in\g_1$ and $y\in \g_0$.
Then 
$x y v'=[xy]v'+yxv'=0$. 
Since $\g_1 v_{\Lambda_i}=0$, the claim follows.
\end{dokaz}

The standard $ \tilde{\mathfrak g}$-module $L(\Lambda_0)$ is a vertex operator algebra and
$L(\Lambda_1),\dots,L(\Lambda_\ell)$ are modules for vertex operator algebar $L(\Lambda_0)$ (cf. \cite{FLM}, \cite{LL}). 
By using Theorem 6.2 in \cite{FF} it is easy to see that the space of intertwining operators
$${\mathcal Y}:L(\Lambda_1) \to \textrm{Hom}\,(L(\Lambda_i),L(\Lambda_{i+1}))\{\{z\}\},
\qquad {\mathcal Y}(w,z)u=\sum_{m\in\mathbb Q}w_m u\,z\sp{-m-1}$$
is $1$-dimensional for $i=0,\dots,\ell-1$.

By Lemma 6.1 in \cite{P4} for such nonzero $\mathcal Y$ there exists $m\in\mathbb Q$ such that
\begin{equation} \label{E: intertwining 1}
V_1\otimes V_i \to V_{i+1}, \qquad w\otimes u\mapsto w_m u
\end{equation}
is a nonzero homomorphism of  $\g_0$-modules. It is easy to see that the multiplicity of $V_{i+1}$ in $V_1\otimes V_i$ is $1$ --- one way to see this is by using Parthasarathy-Ranga Rao-Varadarajan's
theorem 5.2 in \cite{FF} --- hence we can normalize
$\mathcal Y$ so that the map  \eqref{E: intertwining 1} is the homomorphism of  $\g_0$-modules
\begin{equation} \label{E: intertwining 2}
V_1\otimes V_i \to V_{i+1}, \qquad v_m u=v\wedge u.
\end{equation}
From the commutator formula for $\mathcal Y$ and Lemma \ref{Top_rel2} we have:
\begin{prop}
\label{IntOpTop}
For $v\in V_1=\C^\ell$:
\begin{itemize}
\item ${\mathcal Y}(v,z)$ commutes with $\gt_1$
\item for $u\in V_i$, the coefficient of $z\sp{-m-1}$ of ${\mathcal Y}(v,z)u$ is
$$v_m u=v\wedge u.$$
\end{itemize}
\end{prop}

\section{Simple current operator} \label{S: SimpCurr}

Recall that we have fixed the minuscule coweight $\omega=\omega_\ell \in\h^*$.
We shall use simple current operators
$$L(\Lambda_i) \xrightarrow{[\omega]} L(\Lambda_{\ell-i}) \xrightarrow{[\omega]}
L(\Lambda_i) $$
\smallskip
such that simple current commutation relation 
\begin{equation} \label{OmegaComm} x_\alpha(n)[\omega]=[\omega] x_\alpha (n+\alpha(\omega)),
\qquad \alpha\in R, \quad n\in\mathbb Z\end{equation}
holds (see \cite{DLM} and \cite{Ga}, or Remark 5.1 in \cite{P4}). Then
\begin{eqnarray*}
x_\gamma(-n-1) [\omega] & = & [\omega] x_\gamma(-n)\quad
\textrm{for}\quad \gamma \in \Gamma\\
x(\mu)[\omega] & = & [\omega]x(\mu^+) \quad \textrm{for a monomial}\
x(\mu)\in U(\gt_1^-),
\end{eqnarray*}
where by
$x(\mu^+)$ we denote a monomial obtained by  raising degrees in $x(\mu)$ by $1$.
Also,
\begin{eqnarray*}
x_\alpha(-n) [\omega] & = & [\omega] x_\alpha(-n)\quad
\textrm{for}\quad x_\alpha \in \g_0,\\
\sigma [\omega] & = & [\omega] \sigma \quad \textrm{for}\quad \sigma\in\Weyl.
\end{eqnarray*}
\bigskip

In the level $k>1$ case, 
for $\Lambda = k_0\Lambda_0 +
k_1\Lambda_1+ \cdots + k_\ell\Lambda_\ell$,
we embed $L(\Lambda)$ in a tensor
product of standard modules of level $1$
$$L(\Lambda)\subset L(\Lambda_{0})^{k_0}\otimes
L(\Lambda_{1})^{k_1}\otimes \cdots \otimes
L(\Lambda_{\ell})^{k_\ell},$$
with a highest weight vector
$$v_\Lambda=v_{\Lambda_0}^{\otimes k_0}\otimes\dots\otimes
v_{\Lambda_\ell}^{\otimes k_\ell},$$ 
and we use the level $k$ simple current operator $[\omega]\colon L(\Lambda)\to L(\Lambda')$,
$$[\omega]=[\omega]\otimes\dots\otimes [\omega].$$

\section{Relations} \label{S: Rel}

From \eqref{DC_rel} we immediately obtain the following relations: 

\begin{lm} \label{rel4_lm}
On a level one module $L(\Lambda_r)$
\begin{eqnarray}
 \label{rel_iiii} & & \hspace{-3ex} x_{ii}(-1) x_{ii}(-1)  =  0,\\
 \label{rel_iiij} & & \hspace{-3ex} x_{ij}(-1) x_{ii}(-1)   =  0,\quad \textrm{for}\ i<j,\\
 \label{rel_ijjj} & & \hspace{-3ex} x_{jj}(-1) x_{ij}(-1)   =  0,\quad \textrm{for}\ i<j,\\
 \label{rel_iijj} & & \hspace{-3ex} x_{ij}(-1)x_{ij}(-1)  \sim  x_{ii}(-1)x_{jj}(-1), \quad \textrm{for} \ i<j,\\
 \label{rel_iijk} & & \hspace{-3ex} x_{jk}(-1)x_{ii}(-1)  \sim x_{ik}(-1)x_{ij}(-1)  =0, \quad \textrm{for} \ i<j<k,\\
 \label{rel_ijjk} & & \hspace{-3ex} x_{jk}(-1)x_{ij}(-1)  \sim x_{jj}(-1)x_{ik}(-1)  =0, \quad \textrm{for} \ i<j<k,\\
 \label{rel_ijkk} & & \hspace{-3ex} x_{kk}(-1)x_{ij}(-1)  \sim x_{jk}(-1)x_{ik}(-1)  =0, \quad \textrm{for} \ i<j<k,\\
 \label{rel_ijkl} & & \hspace{-3ex} x_{jk}(-1)x_{il}(-1)  \sim C_1 x_{jl}(-1)x_{ik}(-1)  + C_2 x_{kl}(-1)x_{ij}(-1), 
\end{eqnarray}
for $i<j<k<l$ and some $C_1,C_2\in \C^\times$.
\end{lm}

\begin{lm} \label{rel2_lm}
If $i,j\in\{s_1,s_2,\dots,s_r\}$, where $1\leq s_1,\dots,s_r \leq \ell$, then
$$x_{i j}(-1) v_{s_1 s_2 \dots s_r}=0.$$
\end{lm}

\begin{dokaz} Let $i=s_p$, $j=s_q$, $p\leq q\leq r$.
By lemma \ref{IC1_lm}, we have
\begin{equation} \label{eq_r2}
x_{p q}(-1)v_{1 2 \dots r}=0.
\end{equation}
Next, by acting on \eqref{eq_r2} with a linear map $\sigma$ corresponding to $\sigma\in\Weyl$ such that
$\sigma(t)=s_t$, $t=1,\dots,r$ (see \eqref{E: Weyl map}),   we obtain the desired relation.
\end{dokaz}

\begin{lm} \label{rel3_lm}
 If $i\notin\{s_1,s_2,\dots,s_r\}$, then
\begin{equation} \label{eq_r3}
x_{i s_p}(-1) v_{s_1 s_2 \dots s_r} \sim x_{s_p s_p}(-1)v_{s_1 \dots i \dots \underline{s_p} \dots s_r}.
\end{equation}
\end{lm}
 
We use underline to denote that the corresponding indices should be excluded.
 
\begin{dokaz}
Let $s_{q-1}<i<s_q$, $q\leq p$. By lemma \ref{rel2_lm}, we have
\begin{equation} 
x_{q+1, p+1}(-1)v_{1 \dots q-1, q+1 \dots r+1}=0.
\end{equation}
We act with $x_{\alpha_{q}}\in \g_0$ and obtain
\begin{equation} 
x_{q, p+1}(-1)v_{1 \dots q-1, q+1 \dots r+1} - C' x_{q+1, p+1}(-1)v_{1 \dots q, q+2 \dots r+1}=0,
\end{equation}
for some $C'\in \C^\times$. Hence
\begin{equation} 
x_{q, p+1}(-1)v_{1 \dots q-1, q+1 \dots r+1} \sim x_{q+1, p+1}(-1)v_{1 \dots q, q+2 \dots r+1}.
\end{equation}
By induction, we obtain
\begin{equation} \label{eq_r4}
x_{q, p+1}(-1)v_{1 \dots q-1, q+1 \dots r+1}\sim x_{p+1, p+1}(-1)v_{1 \dots p, p+2 \dots r+1}.
\end{equation}
By acting on \eqref{eq_r4} with a linear map $\sigma$ corresponding to $\sigma\in\Weyl$ such that
$\sigma(t)=s_t$, for $t\in\{1,\dots,q-1\}$, $\sigma(t)=s_{t-1}$, for $t\in\{q+1,\dots,r+1\}$, and 
$\sigma(q)=i$, we obtain \eqref{eq_r3}.
\end{dokaz}

Lemma \ref{rel2_lm} can be generalized in the following way:

%

\begin{lm} \label{rel5_lm}
Let $x(\pi)=x_{i_m j_m}(-1) \cdots x_{i_1 j_1}(-1)$, $i_t\leq j_t$, and assume $j_1<\dots< j_m$.
Let $I\subset \{1,\dots,\ell\}$ be such that $\{1,\dots,j_{m}\}\setminus \{j_1,\dots,j_{m}\}\subset I$.
 If $j_m\in I$, then
$$x(\pi) v_I=0.$$
\end{lm}

\begin{dokaz}
If $m=1$ this is just lemma \ref{rel2_lm} since $i_1,j_1\in I$. For $m>1$ we use induction. 
If $i_m\in I$, then lemma \ref{rel2_lm} implies $x(\pi)v_I=0$.
If $i_m\notin I$, then $i_m=j_r$, for some $r<m$. By lemma \ref{rel3_lm}
$$x_{j_r j_m} (-1) v_I= x_{j_m j_m} (-1) v_{I'},$$
where $I'=(I\setminus\{j_m\})\cup\{j_r\}$. By induction,
$$x_{i_r j_r}(-1) \cdots x_{i_1 j_1}(-1) v_{I'} = 0,$$
which gives the claim.
\end{dokaz}

\begin{lm} \label{rel55_lm}
Let $x(\pi)=x_{i_m j_m}(-1) \cdots x_{i_1 j_1}(-1)$, $i_t\leq j_t$, and assume $j_1< \dots<  j_{m-1} \leq j_m$. 
Let $I\subset \{1,\dots,\ell\}$ be such that $\{1,\dots,j_{m}\}\setminus \{j_1,\dots,j_{m}\}\subset I$.
 If $j_{m-1}= j_m$, then
$$x(\pi) v_I=0.$$
\end{lm}

\begin{dokaz}
For $m=2$ we have $x(\pi)v_I = x_{i_2 j_2}(-1) x_{i_1 j_2}(-1) v_I$. Assume $i_1\leq i_2$.
If $i_2=j_2$, then $x(\pi) v_I=0$, by \eqref{rel_ijjj}.
Otherwise, $i_1,i_2\in I$ and $x(\pi) v_I \sim x_{j_2 j_2}(-1) x_{i_1 i_2}(-1) v_I=0$,
by \eqref{rel_ijkk}, \eqref{rel_iijj} and lemma \ref{rel2_lm}.

For $m>2$ we use induction. Assume $i_{m-1}\leq i_m$.
If $i_m=j_m$, then $x(\pi) v_I=0$, by \eqref{rel_ijjj}.
Otherwise, by \eqref{rel_ijkk} and \eqref{rel_iijj},
$$x(\pi)v_I \sim x_{j_m j_m}(-1) x_{i_{m-1} i_m}(-1) x_{i_{m-2} j_{m-2}} \cdots x_{i_1 j_1}(-1) v_I.$$
If $i_m\in I$, then $x(\pi) v_I=0$, by lemma \ref{rel5_lm}.
If $i_m\notin I$, then $i_m=j_r$, for some $r<m-1$. By induction,
$$x_{i_{m-1} j_r}(-1) x_{i_{r} j_{r}} \cdots x_{i_1 j_1}(-1) v_I=0,$$
from which the claim follows.
\end{dokaz}

So far we have described conditions upon which a monomial of ``small'' degree will annihilate vectors from the top of a standard module.
In the following two propositions we describe in more detail the action of these monomials on vectors from the top. 

\begin{prop}\label{rel6_lm}
Let $x(\pi)=x_{i_m j_m}(-1)\cdots x_{i_2 j_2}(-1) x_{i_1 j_1}(-1)$ satisfy
difference conditions, i.e. $i_1<\dots < i_m$, $j_1<\dots <j_m$ and $i_t\leq j_t$.
Set $I=\{1,\dots,\ell\}\setminus \{j_1,\dots,j_m\}$, $I'=\{i_1,\dots,i_m\}$. Then
\begin{equation}
\label{w_rel0} x(\pi) v_I \sim  [\omega] v_{I'}.
\end{equation}
\end{prop}

\begin{dokaz}
We first show 
\begin{equation} \label{w_rel1}
x_{mm}(-1) x_{m-1,m-1}(-1)\dots x_{22}(-1) x_{11}(-1) v_{m+1,\dots,\ell} \sim [\omega] v_{1 2\dots m},
\end{equation}
 i.e. we show that $[\omega]^{-1}x_{mm}(-1) \dots x_{11}(-1) v_{m+1,\dots,\ell}$
is a highest weight vector with the highest weight $\Lambda_m$.
Like in the proof of lemma \ref{IC1_lm}, one easily sees that $x_{mm}(-1) \dots x_{11}(-1) v_{m+1,\dots,\ell}\neq 0$.
By \eqref{OmegaComm} and lemma 
\ref{Top_rel2}
we have
\begin{eqnarray*}
& & \hspace{-13ex}
x_{-\theta}(1)[\omega]^{-1}x_{mm}(-1) \dots x_{11}(-1) v_{m+1,\dots,\ell}\\
& &  = [\omega]^{-1}x_{mm}(-1) \dots x_{11}(-1) x_{-\theta} (2) v_{m+1,\dots,\ell} = 0.
\end{eqnarray*}
By \eqref{Top_vect} and \eqref{OmegaComm} and 
\begin{eqnarray*}
& & \hspace{-13ex}
x_{\alpha_i}(0) [\omega]^{-1}x_{mm}(-1) \dots x_{11}(-1) v_{m+1,\dots,\ell} \\
& &  =  [\omega]^{-1}x_{mm}(-1) \dots x_{11}(-1) x_{\alpha_i}(0)  v_{m+1,\dots,\ell} =0,
\end{eqnarray*}
for $i=m+1,\dots,\ell-1$. 
By \eqref{Top_vect}, \eqref{OmegaComm} and lemma \ref{rel2_lm}
\begin{eqnarray*}
& & \hspace{-13ex} x_{\alpha_m}(0) [\omega]^{-1}x_{mm}(-1) \dots x_{11}(-1) v_{m+1,\dots,\ell} \\
& &  = [\omega]^{-1}x_{mm}(-1) \dots x_{11}(-1) v_{m,m+2,\dots,\ell} = 0,
\end{eqnarray*}
By \eqref{Top_vect}, \eqref{OmegaComm}  and \eqref{rel_iiij}
\begin{eqnarray*}
& & \hspace{-5ex} x_{\alpha_i}(0) [\omega]^{-1}x_{mm}(-1) \dots x_{11}(-1) v_{m+1,\dots,\ell}   =
 [\omega]^{-1}x_{mm}(-1) \dots x_{ii}(-1)\\
 & &
 \cdot x_{i+2,i+2}(-1) \dots x_{11}(-1)
(x_{i+1,i+1}(-1) x_{\alpha_i}(0) + x_{i,i+1}(-1))v_{m+1,\dots,\ell} = 0,
\end{eqnarray*}
for $i=1,\dots,m-1$. 
Since $\alpha_\ell=(\ell \ell)$,  by  \eqref{OmegaComm} and lemma \ref{rel2_lm}
\begin{eqnarray*}
& & \hspace{-13ex}
 x_{\alpha_\ell}(0) [\omega]^{-1}x_{mm}(-1) \dots x_{11}(-1) v_{m+1,\dots,\ell}  \\
 & &  
= [\omega]^{-1} x_{mm}(-1) \dots x_{11}(-1) x_{\ell \ell}(-1)  v_{m+1,\dots,\ell} = 0.
\end{eqnarray*}
Hence, \eqref{w_rel1} holds. 

Next, $\Weyl$-action on \eqref{w_rel1} gives
\begin{equation}\label{w_rel2}
x_{j_m j_m}(-1)\cdots x_{j_2 j_2}(-1) x_{j_1 j_1}(-1)v_I \sim [\omega] v_{j_1\dots j_m}.
\end{equation}

Finally, $\g_0$-action on \eqref{w_rel2} gives the claim. Assume that we have shown
\begin{eqnarray}\label{w_rel3}
& &  \hspace{-3ex}
x_{j_m j_m}(-1)\cdots x_{j_{t+1} j_{t+1}}(-1) x_{n j_{t}}(-1) x_{i_{t-1} j_{t-1}}(-1) \cdots x_{i_1 j_1}(-1)v_I \\
\nonumber & &  \hspace{10ex} \sim [\omega] v_{i_1\dots i_{t-1} n j_{t+1}\dots j_m},
\end{eqnarray}
where $i_t<n \leq j_t$. 
Since $[x_{\alpha_{n-1}} (0),x_{n j_{t}}(-1)]=x_{n-1, j_{t}}(-1)$, we act on \eqref{w_rel3} with $x_{\alpha_{n-1}} (0)$. 
We claim that what we will get is
\begin{eqnarray}\label{w_rel4}
& &  \hspace{-3ex}
x_{j_m j_m}(-1)\cdots x_{j_{t+1} j_{t+1}}(-1) x_{n-1, j_{t}}(-1) x_{i_{t-1} j_{t-1}}(-1) \cdots x_{i_1 j_1}(-1)v_I \\
\nonumber & &  \hspace{10ex} \sim [\omega] v_{i_1\dots i_{t-1}, n-1, j_{t+1}\dots j_m}.
\end{eqnarray}
Note that $[x_{\alpha_{n-1}} (0),x_{j_r j_{r}}(-1)]=0$, for  $t<r\leq m$. 

First, consider the case when $n\in I$. In this case also $[x_{\alpha_{n-1}} (0),x_{i_r j_{r}}(-1)]=0$, for $1\leq r< t$. 
If $n-1\in I$, then $x_{\alpha_{n-1}} (0) v_I=0$ and \eqref{w_rel4} follows.
If $n-1\notin I$, then $x_{\alpha_{n-1}} (0) v_I=v_{I'}$, where $I'=(I\setminus \{n\})\cup \{n-1\}$.
Since $n-1\notin I$, then $n-1=j_s$, for some $s<t$. By lemma \ref{rel5_lm}, $$x_{i_{s} j_{s}}(-1) \cdots x_{i_1 j_1}(-1)v_{I'}=0,$$
and \eqref{w_rel4} follows.

Now, consider the case when $n\notin I$. In this case $n=j_s$, for some $r<t$. Then $x_{\alpha_{n-1}} (0) v_I=0$,
 $[x_{\alpha_{n-1}} (0),x_{i_r j_{r}}(-1)]=0$, for $1\leq r< t, r\neq s$, and  $[x_{\alpha_{n-1}} (0),x_{i_s j_{s}}(-1)]=x_{i_s, n-1}(-1)$.
If $n-1\in I$, then, by lemma \ref{rel5_lm},
$$x_{i_{s}, n-1}(-1) x_{i_{s-1} j_{s-1}}(-1) \cdots x_{i_1 j_1}(-1)v_{I}=0.$$
If $n-1\notin I$, then $n-1=j_{s-1}$ and  
$$x_{i_{s} j_{s-1}}(-1) x_{i_{s-1} j_{s-1}}(-1) \cdots x_{i_1 j_1}(-1)v_{I}=0,$$ 
by lemma \ref{rel55_lm}.
In both cases \eqref{w_rel4} follows.

After finitely many steps we will reach \eqref{w_rel0}
\end{dokaz}

\begin{prop}\label{rel1_lm}
Let $x(\pi)=x_{i_m j_m}(-1)\cdots x_{i_2 j_2}(-1) x_{i_1 j_1}(-1)$ be such
that $j_1<\dots < j_m$, $i_r\neq i_{t}$ for $r\neq t$, and $i_t\leq j_t$.
Set $I=\{1,\dots,\ell\}\setminus \{j_1,\dots,j_m\}$, $I'=\{i_1,\dots,i_m\}$. Then
$x(\pi) v_I \sim  [\omega] v_{I'}$.
\end{prop}

\begin{dokaz}
Since $x(\pi)$ does not satisfy difference conditions, we can use relations \eqref{DC_rel}
to express $x(\pi)$ in terms of monomials that satisfy difference conditions. 
We show that in each step of the reduction of $x(\pi)$ we will obtain another monomial $x(\pi')$ such that
$x(\pi)<x(\pi')$, $x(\pi)v_I \sim x(\pi')v_I$ and such that colors of factors of $x(\pi')$ again lie in 
rows $j_1,\dots,j_m$ and columns $i_1,\dots,i_m$, as are the colors of $x(\pi)$. In the end
we will obtain $x(\pi)v_I=x_{i_{s_m} j_m}(-1)\cdots x_{i_{s_1} j_1}(-1)v_I$, where 
$i_{s_1}<\dots < i_{s_m}$, and proposition \ref{rel6_lm} gives the claim.

Let $x(\pi)=x(\pi_1)x(\pi_2)$, where $x(\pi_1)$ is a leading term of some relation \eqref{DC_rel}.
In fact, since colors of factors of $x(\pi)$ all lie in different rows and different columns, relations that
we use are \eqref{rel_ijjk} and \eqref{rel_ijkl}.

First, assume that the reduction is made upon relation \eqref{rel_ijjk}. In this case
$x(\pi_1)v_I \sim  x(\pi_1')v_I$, with $x(\pi_1)<x(\pi_1')$, and after the reduction we obtain another monomial
$x(\pi')=x(\pi_1')x(\pi_2)$ whose colors lie in the same rows and columns as colors of $x(\pi)$ and 
$x(\pi_1)<x(\pi_1')$.

Next, assume that the reduction is made upon relation \eqref{rel_ijkl}. 
Let $x(\pi_1)=x_{i_t j_t}(-1) x_{i_r j_r}(-1)$,  for some $r<t$ and $i_t< i_r < j_r < j_t$. Let 
$$x(\pi_1')=x_{i_r j_t}(-1)x_{i_t j_r}(-1), \quad  
x(\pi_1'')=x_{i_t i_r}(-1)x_{j_r j_t}(-1).$$ 
Then, by \eqref{rel_ijkl}, 
$$x(\pi_1)v_I=C_1 x(\pi_1')v_I+ C_2 x(\pi_1'')v_I,$$
 for some $C_1,C_2\in\C^\times$ and $x(\pi_1)<x(\pi_1'),x(\pi_1'')$. 
Denote by $x(\pi')=x(\pi_1')x(\pi_2)$ and $x(\pi'')=x(\pi_1'')x(\pi_2)$. Then 
$$x(\pi)v_I=C_1 x(\pi')v_I+ C_2 x(\pi'')v_I,$$ and $x(\pi)<x(\pi'),x(\pi'')$.
Note that colors of $x(\pi')$ lie in the same rows and columns as colors of $x(\pi)$.
We claim that $x(\pi'')v_I=0$. Let $n$ be such that $j_{n}\leq i_r < j_{n+1}$.
If $i_r\in I$, then 
$$x_{i_t i_r}(-1)  x_{i_n j_n}(-1)\cdots  x_{i_1 j_1}(-1)v_I=0,$$
 by lemma \ref{rel5_lm}. 
If $i_r\notin I$, then $j_{n}= i_r$ and 
$$x_{i_t i_r}(-1)  x_{i_n j_n}(-1)\cdots  x_{i_1 j_1}(-1)v_I=0,$$ by lemma \ref{rel55_lm}. 
In  both cases we get $x(\pi'')v_I=0$.
\end{dokaz}

\begin{napomena}\label{Switching_Rem} \emph{
Let $x_{i_m j_m}(-1) \cdots x_{i_1 j_1}(-1)$ satisfy difference conditions.
Assume  
\begin{equation} \label{r5_eq}
x_{i_m j_m}(-1) \cdots x_{i_1 j_1}(-1) v_{s_1 \dots s_{\ell-m}}\neq 0.
\end{equation}
If $\{j_1,\dots,j_m\}\cap \{s_1,\dots,s_{\ell-m}\}=\emptyset$, then,
by proposition \ref{rel6_lm}, $$x_{i_m j_m}(-1) \cdots x_{i_1 j_1}(-1) v_{s_1 \dots s_{\ell-m}}
=[\omega] v_{i_1\dots i_m}.$$
If $\{j_1,\dots,j_m\}\cap \{s_1,\dots,s_{\ell-m}\}\neq \emptyset$, we can use lemma \ref{rel3_lm} 
to ``switch'' indices and reduce the expression \eqref{r5_eq} to the
one appearing in proposition \ref{rel1_lm}. 
Concretely: let $t$ be the smallest possible so that
$j_t \in \{s_1,\dots,s_{\ell-m}\}$. Then $i_t\notin \{s_1,\dots,s_{\ell-m}\}$, by \eqref{r5_eq} and lemma \ref{rel2_lm},  and
$$x_{i_t j_t}(-1) v_{s_1 s_2 \dots s_{\ell-m}} \sim  x_{j_t j_t}(-1) 
v_{s_1 \dots i_t \dots \underline{j_t} \dots s_{\ell-m}},$$
by lemma \ref{rel3_lm}.
Proceed inductively -- reduce the expression
$$x_{i_m j_m}(-1) \cdots x_{j_t j_t}(-1) \cdots x_{i_1 j_1}(-1) 
v_{s_1 \dots i_t \dots \underline{j_t} \dots s_{\ell-m}}$$ 
to the one appearing in proposition \ref{rel1_lm}.\\
\indent Hence, in this case, there is no index occuring more than twice in the sequence 
$i_1,\dots,i_m$,\, $j_1,\dots,j_m$,\, $s_1,\dots,s_{\ell-m}$, and 
$$ 
x_{i_m j_m}(-1) \cdots x_{i_1 j_1}(-1) v_{s_1 \dots s_{\ell-m}} \sim [\omega]v_{r_1 \dots r_m},
$$
where $r_1,\dots,r_m$ are exactly those indices appearing twice in the sequence.
}\end{napomena}

\section{Proof of linear independence} \label{S: LinIndep}

Let $x(\pi)=x_{i_m j_m}(-1)\cdots x_{i_2 j_2}(-1) x_{i_1 j_1}(-1)$ satisfy difference and initial conditions
on $W(\Lambda_i)$. Set $J=\{1,\dots,\ell\}\setminus \{j_1,\dots, j_m \}$, $I=\{i_1,\dots,i_m\}$.
By Proposition \ref{IntOpTop}
there are 
operators, denoted by $w_1, w_2$ that
commute with $\gt_1$ and such that
\begin{equation} \label{DefOpIsp_eq}
w_1 v_{\Lambda_i}=v_J,\quad w_2 v_I= v_{\Lambda_{i_m}}.
\end{equation}
Moreover, operators $w_1$ and $w_2$ act on $V_{i}$ and $V_m$, correspondingly, as multiplication
in the exterior algebra by suitable vectors. Let $w_2^{\omega} = [\omega]w_2 [\omega]^{-1} $; it also commutes with $\gt_1$.
Then, by proposition \ref{rel6_lm}, 
$$w_2^\omega w_1 x(\pi) v_{\Lambda_i} \sim [\omega] v_{\Lambda_{i_m}}.$$

Let $x(\mu)=x_{r_{n} s_{n}}(-1)\cdots x_{r_2 s_2}(-1) x_{r_1 s_1}(-1)$ also satisfy difference and initial conditions
on $W(\Lambda_i)$, and let $x(\mu)>x(\pi)$. We will  show that
\begin{equation}\label{Annih_eq}
w_2^\omega  w_1 x(\mu) v_{\Lambda_i}=0.
\end{equation} 

If $n>m$, then either $$x_{r_{m} s_{m}}(-1)\cdots x_{r_1 s_1}(-1) v_J=0$$ 
or
$$x_{r_{m} s_{m}}(-1)\cdots x_{r_1 s_1}(-1) v_J \sim  [\omega] v',$$ for some $v'\in V_m$
(cf. Remark \ref{Switching_Rem}). In the first case relation \eqref{Annih_eq} directly follows.
In the second case we have
$$ x(\mu)v_J \sim [\omega]x_{r_{n} s_{n}}(0)\cdots x_{r_{m+1} s_{m+1}}(0) v'=0,$$
by lemma \ref{Top_rel2}.
Hence, \eqref{Annih_eq} follows.

If $n=m$, then there exists $t$, $1\leq t \leq m$, such that
$(r_1 s_1)=(i_1 j_1)$, $\dots$, $(r_{t-1} s_{t-1})= (i_{t-1} j_{t-1})$
and $(r_{t} s_{t})> (i_{t} j_{t})$. Then either $r_t=i_t$ and $s_t<j_t$,
or $r_t < i_t$.

Consider first the case when $r_t=i_t$ and $s_t<j_t$. Then $j_{t-1}=s_{t-1}<s_t<j_t$, hence
$s_t\in J$. Then $x_{r_t s_t}(-1) \cdots x_{r_1 s_1}(-1) v_J=0$, by lemma \ref{rel5_lm}, 
and \eqref{Annih_eq} follows.

Now, consider the case when $r_t < i_t$. If $x(\mu)v_J=0$, we are done. 
If $x(\mu)v_J\neq 0$, then, by remark \ref{Switching_Rem},  in the multiset 
$\{r_1,\dots,r_m\}\cup \{s_1,\dots,s_m\}\cup J$ there is no index occuring more than twice.
Furthermore, since $r_t<j_t$ and $s_p=j_p,$ for $p<t$, index $r_t$ appears twice in the forementioned multiset.
Hence $x(\mu)v_J \sim  [\omega] v_{I'}$ and $r_t\in I'$. Since $r_t\notin I$ and $r_t < i_t$, the action of $w_2$ will annihilate $v_{I'}$ (see \eqref{DefOpIsp_eq} and proposition \ref{IntOpTop}), i.e. $w_2^\omega [\omega] v_{I'}=0$.

Set $w=w_2^\omega w_1$. By the above considerations we have

\begin{prop} \label{OpIsp1_prop}
Let  $x(\pi)$ satisfy difference and initial
conditions for $L(\Lambda_i)$. Write
$x(\pi)=x(\pi_1)x(\pi_2)$, where $x(\pi_1)$ is the $(-1)$-part of a
monomial, and $x(\pi_2)$ the rest of the monomial. Then there
exists an operator $w:L(\Lambda_i) \to L(\Lambda_{i'})$ such that:
\begin{itemize}
\item $w$ commutes with $\gt_1$,
\item $w x(\pi_1) v_{\Lambda_i}  \sim  [\omega] v_{\Lambda_{i'}},$
\item $x(\pi_2^+)$ satisfies IC and DC for $L(\Lambda_{i'})$,
\item if  $x(\pi')$ has a $(-1)$-part $x(\pi_1')$  greater than $x(\pi_1)$, then $w x(\pi') v_{\Lambda_i}=0$.
\end{itemize}
\end{prop}

\begin{dokaz}
We have already shown the majority of the claims. It remains to see that $x(\pi_2^+)$ satisfies IC and DC for $L(\Lambda_{i'})$, but this is clear from
the description of level $1$ difference conditions in \eqref{DClev1_eq} and the definition of initial conditions (see also remark \ref{ICred_nap}). 
\end{dokaz}	

Like in the $A_\ell^{(1)}$-case, proposition \ref{DCFact_prop} enables us to straightforwardly generalize
proposition \ref{OpIsp1_prop} for higher levels (cf. [T2]).
Also, the proof of linear independence is the same as in the $A_\ell^{(1)}$-case, hence we give here
only the sketch of the proof.


\begin{tm}\label{FSbaza_tm}
The set
$$\{x(\pi)v_\Lambda \,|\, x(\pi) \ \textrm{satisfies DC and IC for}\ L(\Lambda)\}$$
is a basis of $W(\Lambda)$.
\end{tm}

\begin{sketch}
Assume
\begin{equation}\label{LinZav_jed}\sum c_\pi x(\pi)v_\Lambda=0,\end{equation}
where all monomials $x(\pi)$ satisfy difference and initial
conditions for $W(\Lambda)$. Fix $x(\pi)$ in \eqref{LinZav_jed} and assume that
$c_{\pi'}=0$  for $x(\pi')<x(\pi)$.
We show that $c_\pi=0$.

Choose an operator $w$ from proposition \ref{OpIsp1_prop} and apply it
on \eqref{LinZav_jed}. By proposition \ref{OpIsp1_prop}, we get  
\begin{eqnarray*}
0 & = & w \sum c_{\pi'}x(\pi') v_\Lambda\\ & = & 
w \sum_{\pi_1'>\pi_1}c_{\pi'}x(\pi') v_\Lambda +w \sum_{\pi_1'<\pi_1}c_{\pi'}x(\pi')
 v_\Lambda+ w \sum_{\pi_1'=\pi_1}c_{\pi'}x(\pi') v_\Lambda\\
 & = & w \sum_{\pi_1'=\pi_1}c_{\pi'}x(\pi') v_\Lambda 
 \sim \sum_{\pi_1'=\pi_1}c_{\pi'}x(\pi_2') [\omega] v_{\Lambda'}\\ & = & 
 [\omega] \sum_{\pi_1'=\pi_1}c_{\pi'}x(\pi_2'^+) v_{\Lambda'}
\end{eqnarray*}
Since $[\omega]$ is injective, it follows that
$$ \sum_{\pi_1'=\pi_1}c_{\pi'}x(\pi_2'^+)
v_{\Lambda'}=0.$$  This is a relation of linear dependence 
between monomial vectors satisfying difference and initial conditions for  $W(\Lambda')$
with all monomials of degree greater than the degree of $x(\pi)$. By the induction hypothesis 
they are linearly independent, and, in particular, $c_\pi=0$.
\end{sketch}

\end{document}